\documentclass{amsart}%{proc-l}

\theoremstyle{remark}
\newcommand{\ip}[2]{\mbox{$(#1|#2)$}}
\newcommand{\CC}{{\bf C}}
\newtheorem{lemma}{Lemma}[section]
\newtheorem{theorem}{Theorem}%[section]

\newtheorem{remark}[lemma]{Remark}
\newtheorem{corollary}[lemma]{Corollary}
\newtheorem{problem}{Problem}

\newcommand{\csa}{$C^*$-algebra}

\newcommand{\jcst}{$JC^*$-triple}
\newcommand{\jwst}{$JW^*$-triple}
\begin{document}

\title[Representation of %contractively complemented
Hilbertian operator spaces on the Fock space]{Representation of
contractively complemented Hilbertian operator spaces on the Fock
space}

%    Information for first author
\author{Matthew Neal}
%    Address of record for the research reported here
\address{Department of Mathematics, Denison University, Granville,
Ohio 43023}
%    Current address
%\curraddr{Department of Mathematics, Denison University,
%Granville, Ohio 43023}
 \email{nealm@denison.edu}
%    \thanks will become a 1st page footnote.
%\thanks{The first author was supported in part by .}

%    Information for second author
\author{Bernard Russo}
\address{Department of Mathematics, University of California,
Irvine, California 92697-3875}
 \email{brusso@math.uci.edu}
%\thanks{Support information for the second author.}

%    General info

\date{January 1, 1994 and, in revised form, June 22, 1994.}

%\dedicatory{This paper is dedicated to our authors.}


\begin{abstract}
The operator spaces $H_n^k$ $1\le k\le n$, generalizing the row
and column Hilbert spaces, and arising in the authors' previous
study of contractively complemented subspaces of $C^*$-algebras,
are shown to be homogeneous and completely isometric to a space of
creation operators on a subspace of the anti-symmetric Fock space.
The completely bounded Banach-Mazur distance from $H_n^k$ to row
or column space is explicitly calculated.
\end{abstract}

\maketitle
\section*{Introduction and Preliminaries}
A  well-known result of Friedman and Russo (\cite[Theorem
2]{FriRus85}) states that if a subspace $X$ of a \csa\ $A$  is the
range of a contractive projection on $A$, then $X$ is isometric to
a \jcst, that is, a norm closed subspace of $B(H,K)$ stable under
the triple product $ab^*c+cb^*a$. If $X$ is atomic (in particular,
finite-dimensional), then it is isometric to a direct sum of
Cartan factors of types 1 to 4.

The authors showed in \cite{NeaRus03} that this latter result
fails, as it stands, in the category of operator spaces. In that
paper, we defined a family of $n$-dimensional Hilbertian operator
spaces $H_n^k$, $1\le k\le n$, generalizing the row and column
Hilbert spaces $R_n$ and $C_n$ and showed that in the above
result, if $X$ is atomic, the word ``isometric'' can  be replaced
by ``completely semi-isometric,'' provided the spaces $H_n^k$ are
allowed as summands along with the Cartan factors (\cite[Theorem
2]{NeaRus03}). It is pointed out in \cite{NeaRus03} that  the
space $H_n^k$ is contractively complemented in some $B(K)$, and
for $1<k<n$, is not completely (semi-)isometric to either of the
Cartan factors $B(\CC,\CC^n)=H_n^1$ or $B(\CC^n,\CC)=H_n^n$, and
that these spaces appeared in a slightly different form and
context in \cite{AraFri78}. It is also shown in \cite[Theorem
3]{NeaRus03} that finite dimensional \jcst s which are
contractively complemented in a \csa\ can be classified up to
complete isometry.

In this paper, we study the operator space structure of the spaces
$H_n^k$. Besides being a generalization of the row and column
Hilbert spaces, as shown in Lemma~\ref{lem:2.1} below, they are
completely isometric to the span of creation operators on a
subspace of the anti-symmetric Fock space. Thus they are related
to the operator space denoted by $\Phi_n$ in \cite[section
9.3]{Pisier03}, which is the span of the creation operators on the
full anti-symmetric Fock space. $\Phi_n$ is the unique operator
space which is completely isometric to the span of $n$ operators
satisfying the canonical anticommutation relations (CAR),
\cite[Theorem 9.3.1]{Pisier03}, and $\cap_{k=1}^n H_n^k$ is
completely isometric to $\Phi_n$. We show in
 Theorem~\ref{thm:2} below that all finite dimensional Hilbertian operator spaces $X$ which are
contractively complemented in some \csa\ are completely isometric
to the diagonal of two spaces, one space being an intersection of
some of
 the spaces $H_n^k$ and the
other space lying in the kernel of the projection which maps onto
$X$. Since any intersection of the spaces $H_n^k$ is also
completely isometric to a space of creation operators on a
subspace of the full anti-symmetric Fock space,
Theorem~\ref{thm:2} can be interpreted as saying that every
contractively complemented Hilbertian operator space is, up to
complete isometry,
 essentially a space of creation operators. This result is analogous to the result of
 Robertson, \cite{Robertson91}, which states that every {\it completely} contractively complemented Hilbertian operator
 space is completely isometric to either row or column space.

The operator space structures of the row and column Hilbert spaces
$R_n$ and $C_n$ have been well studied, and in particular it is
known that they are homogeneous, dual to each other in the
operator space sense, and have completely bounded Banach-Mazur
distance $n$ between them.  We show here
 that $H_n^k$ is homogeneous (Theorem~\ref{thm:1})
 and we give an explicit formula for the completely bounded Banach-Mazur distance
from it to $R_n=H_n^n$ and $C_n=H_n^1$ (Theorem~\ref{thm:3}). This
answers a question we posed in \cite{NeaRus03} and shows,
interestingly, that the points $R_n$, $C_n$ and $H_n^k$ lie on a
straight line in the metric space of all operator spaces of
dimension $n$.

Recall that a Cartan factor of type 1 is $B(H,K)$ for complex
Hilbert spaces $H$ and $K$. To define the Cartan factors of types
2 and 3, fix a conjugation $J$ on a complex Hilbert space $H$,
that is, a conjugate-linear isometry of order 2, and for $x\in
B(H)$, let $x^t=Jx^*J$. A Cartan factor of type 2 (respectively of
type 3) is
 $A(H,J)
=\{x\in B(H):x^t=-x\}$ (respectively $S(H,J)=\{x\in
B(H):x^t=-x\}$). A Cartan factor of type 4 is the spin factor (cf.
\cite[Subsection 3.1]{NeaRus03}).

An {\it operator space} is a subspace $X$ of $B(H)$, the space of
bounded linear operators on a complex Hilbert space. Its {\it
operator space structure} is given by the sequence of norms on the
set of matrices $M_n(X)$ with entries from $X$, determined by the
identification $M_n(X) \subset M_n(B(H))=B(H\oplus H \oplus \cdots
\oplus H)$. See \cite{Pisier03} for the general theory of operator
spaces, which is now extensive and covered in several other
monographs, for example \cite{EffRau00}, \cite{Paulson03}, and the
forthcoming \cite{BleLeM04}.  Let us just recall that a linear
mapping $\varphi:X\rightarrow Y$ between two operator spaces is
{\it completely bounded} if the induced mappings
$\varphi_n:M_n(X)\rightarrow M_n(Y)$ defined by
$\varphi_n([x_{ij}])=[\varphi(x_{ij})]$ satisfy
$\|\varphi\|_{\mbox{cb}}:=\sup_n\|\varphi_n\|<\infty$. A
completely bounded map is a {\it completely bounded isomorphism}
if its inverse exists and is completely bounded. Two operator
spaces are {\it completely isometric} if there is a linear
isomorphism $T$ between them with $\|T\|_{\mbox{cb}}=
\|T^{-1}\|_{\mbox{cb}}=1$. We call $T$ a {\it complete isometry}
in this case.

In the matrix representation for $B(\ell^2)$ consider the {\it
column Hilbert space} $C=\overline{\mbox{sp}}\{e_{i1}:i\ge 1\}$
and the {\it row Hilbert space}
$R=\overline{\mbox{sp}}\{e_{1j}:j\ge 1\}$ and their finite
dimensional versions $C_n=\mbox{sp}\{e_{i1}:1\le i\le n\}$ and
$R_n=\mbox{sp}\{e_{1j}:1\le j\le n\}$. Here of course $e_{ij}$ is
the operator defined by the matrix with a 1 in the $(i,j)$-entry
and zeros elsewhere. Although $R$ and $C$ are Banach isometric,
they are not completely isomorphic; and $R_n$ and $C_n$, while
completely isomorphic,
 are not completely isometric.

  An operator space is said to be
{\it homogeneous} if every bounded linear map on it is completely
bounded with the norm and completely bounded norm coinciding (see
\cite[9.2]{Pisier03}) and it is {\it Hilbertian} if it is
isometric to a Hilbert space. A linear map of one operator space
into another is said to be a {\it complete semi-isometry} if it is
isometric and completely contractive. The {\it completely bounded
Banach-Mazur distance} between two (completely isomorphic)
operator spaces $E,F$ is defined by
\[
\mbox{d}_{\mbox{cb}}(E,F)=\inf\{\|u\|_{\mbox{cb}}
\|u^{-1}\|_{\mbox{cb}}:u:E\rightarrow F\mbox{ complete isomorphism
}\}.
\]

This paper is organized as follows. In section~\ref{sec:1}, we
show that the spaces $H_n^k$ are homogeneous operator spaces.
Although we  use some multilinear algebra, our proof is direct and
does not make use of the identification of $H_n^k$ with a space of
creation operators. In section~\ref{sec:22} we establish the
complete isometry of $H_n^k$ with a space of creation operators
and use it to describe the fine structure of the range of a
contractive projection on a \csa\ in case said range
 is isometric to a Hilbert space. We also
establish some spectral properties of creation operators. In
section~\ref{sec:2} we compute explicitly the completely bounded
Banach-Mazur distance from the space $H_n^k$ to the column and row
Hilbert spaces $H_n^1$ and $H_n^n$ and state some problems for
further study.

%% The correct journal style for \specialsection is all uppercase; a known bug
%% in amsart.cls prevents this, so input must be uppercase until it is fixed.
%\specialsection*{This is a Special Section Head}
%\specialsection*{THIS IS A SPECIAL SECTION HEAD}
%This is an example of a special section head%
%%%%%%%%%%%%%%%%%%%%%%%%%%%%%%%%%%%%%%%%%%%%%%%%%%%%%%%%%%%%%%%%%%%%%%%%
%\footnote{Here is an example of a footnote. Notice that this footnote
%text is running on so that it can stand as an example of how a footnote
%with separate paragraphs should be written.
%\par
%And here is the beginning of the second paragraph.}%
%%%%%%%%%%%%%%%%%%%%%%%%%%%%%%%%%%%%%%%%%%%%%%%%%%%%%%%%%%%%%%%%%%%%%%%%

\section{Homogeneity of the spaces $H_n^k$}\label{sec:1}%\section{This is a numbered first-level section head}
We begin by recalling from \cite[Sections 6,7]{NeaRus03} the
construction of the spaces $H_n^k$, $1\le k\le n$. Let $I$ denote
a subset of $\{1,2,\ldots,n\}$ of cardinality $|I|=k-1$. The
number of such $I$ is $q:={n\choose k-1}$. Let $J$ denote a subset
of $\{1,2,\ldots,n\}$ of cardinality $|J|=n-k$. The number of such
$J$ is $p:={n\choose n-k}$. Unless otherwise noted, we shall
assume that each $I=\{i_1,\ldots,i_{k-1}\}$ is such that
$i_1<\cdots <i_{k-1}$, and that the collection
$\{I_1,\ldots,I_q\}$ of all such subsets is ordered
lexicographically. Similarly, if $J=\{j_1,\ldots,j_{n-k}\}$, then
$j_1<\cdots<j_{n-k}$ and $\{J_1,\ldots,J_p\}$ is ordered
lexicographically.

We shall use the notation $e_i$ to denote the column vector with a
1 in the $i^{\mbox{th}}$ position and zeros elsewhere. Thus
$e_1,\ldots,e_n$ denotes the canonical basis of column vectors for
$\CC^n$, and for example $e_{J_1},\ldots,e_{J_p}$ denotes the
canonical basis of column vectors for
$\CC^p$. %The vector $e_J$ will also denote the wedge product $e_{j_1}\wedge\cdots\wedge
%e_{j_{n-k}}$, and with this convention, $e_{J_1},\ldots,e_{J_p}$ is an orthonormal
%basis for $\wedge^{n-k}\CC^n$. The meaning of $e_J$ will be clear from its context.

The space $H_n^k$ is the linear span of matrices $b_i^{n,k}$,
$1\le i\le n$, given by
\[
b_i^{n,k}=\sum_{I\cap J=\emptyset,(I\cup
J)^c=\{i\}}\epsilon(I,i,J)e_{J,I},
\]
where $e_{J,I}=e_J\otimes e_I=e_Je_I^t\in
M_{p,q}(\CC)=B(\CC^q,\CC^p)$, and $\epsilon(I,i,J)$  is the
signature of the permutation taking
$(i_1,\ldots,i_{k-1},i,j_1,\ldots,j_{n-k})$ to
$(1,\ldots,n)$.\footnote{In \cite{NeaRus03}, $\epsilon(I,i,J)$ is
also denoted by $\epsilon(I,J)$. However, in this paper,
$\epsilon(I,J)$ will denote the signature of the permutation
taking $(i_1,\ldots,i_{k-1},j_1,\ldots,j_{n-k})$ to
$(1,\ldots,\hat{i},\ldots,n)$} Since the $b_i^{n,k}$ are the image
under a triple isomorphism (actually ternary isomorphism) of a
rectangular grid in a \jwst\ of rank one, they form an orthonormal
basis for $H_n^k$ (cf. the beginning of subsection 5.3 and the
beginning of section 7 of \cite{NeaRus03}).

In the rest of this section, we shall use the following lemma
about determinants, whose proof can be found, for example, in
\cite{Yokonuma92}.

\begin{lemma}\label{lem:det}
Let $X=[\xi_{ij}]$ be an $n\times m$ matrix. Let
$H\subset\{1,\ldots,n\}$ and $K\subset\{ 1,\ldots,m\}$ both have
cardinality $r\le\min\{n,m\}$. Let $X_{H,K}$ denote the
corresponding $r\times r$ submatrix.
\begin{description}
\item[(i)] If $x_i=\sum_{j=1}^n\xi_{ji}e_j\in\CC^n$ for $1\le p
\le n$, then $x_1\wedge\cdots\wedge x_p=\sum_H\det X_{H,L}e_H$,
where the sum is over all $H$ of cardinality $p$,
$L=\{1,\ldots,p\}$, and $X$ is the $n\times p$ matrix
$[\xi_{ij}]$. {\rm (Prop. 3.3, page 84 of \cite{Yokonuma92})}
\item[(ii)] If $X$ is an $n\times n$ matrix and $H\subset \{1,\ldots,n\}$, let $H'$ denote
 the complement of $H$. Then
$\det X=\epsilon(H,H')\sum_R\epsilon(R,R')\det X_{R,H}\det
X_{R',H'}$, where the sum is over all sets $R$ having the same
cardinality of $H$. {\rm (Prop. 3.4(1) page 87 of
\cite{Yokonuma92})}
\item[(iii)] If $H,K\subset\{1,\ldots,n\}$ have cardinality $r$ and $n-r$, and $H\cap K\not=\emptyset$,
then $\sum_R\epsilon(R,R')\det X_{R,H}\det X_{R',K}=0$, where the
sum is over all sets $R$ having  cardinality $r$. {\rm (Prop.
3.4(2), page 87 of \cite{Yokonuma92})}
\end{description}
\end{lemma}

Let $e_1,\ldots,e_n$ be the canonical basis for the column Hilbert
space $C_n=M_{n,1}(\CC) =B(\CC,\CC^n)$ and define an isometry
$\psi:C_n\rightarrow H_n^k$ via $\psi(e_i)=b_i^{n,k}$, $1\le i\le
n$. Let $u=[\begin{array}{ccc}u_1&\cdots&u_n\end{array}]$ be a
unitary matrix so that $u_1,\ldots,u_n$ is an orthonormal basis
for $C_n$. Then, with $u_i=\sum_{j=1}^nu_{ji}e_j$, we have
\[
u=\left[\begin{array}{lll}
u_{11}&\cdots&u_{1n}\\
\vdots&\cdots&\vdots\\
u_{n1}&\cdots&u_{nn}
\end{array}\right],
\]
and
\[
\psi(u_i)=\sum_{j=1}^nu_{ji}b_j^{n,k}=\sum_{j=1}^nu_{ji}\sum_{I\cap
J=\emptyset,(I\cup J)^c=\{i\}}\epsilon(I,i,J)e_{J,I}.
\]

\medskip

\begin{lemma}\label{lem:1.1}
The $(J',I')$-entry of the $p\times q$ matrix $\psi(u_i)$ is given
by
\begin{equation}\label{eq:399}
(\psi(u_i))_{J',I'}=\sum_{I\cap J=\emptyset,(I\cup
J)^c=\{i\}}\epsilon(I,i,J)\frac{\det \overline{u}_{J',J}\det
\overline{u}_{I',I}} {\det \overline{u}},
\end{equation}
where $\overline{u}$ is the complex conjugate of $u$.
\end{lemma}
\proof\ Let us first calculate the left  side of (\ref{eq:399}):
\begin{eqnarray*}
(\psi(u_i))_{J',I'}&=&e_{J'}^t\psi(u_i)e_{I'}\\
&=&\sum_{j=1}^nu_{ji}\sum_{I\cap J=\emptyset,(I\cup
J)^c=\{j\}}\epsilon(I,j,J)e_{J'}^te_Je_I^te_{I'}\\
&=&\left\{ \begin{array}{ll}
0&J'\cap I'\not =\emptyset\\
u_{li}\epsilon(I',l,J')& J'\cap I'=\emptyset,\ (I'\cup
J')^c=\{l\}. \end{array}\right.
\end{eqnarray*}
Before calculating the right  side of (\ref{eq:399}), note that
$\epsilon(I,i,J)\epsilon(I,J)=(-1)^{i+k}$; indeed,
\begin{eqnarray*}
\epsilon(I,i,J)\epsilon(I,J)&=&
(-1)^{k-1}\epsilon(i,I,J)\epsilon(I,J)\\
&=&(-1)^{k-1}\epsilon(I,J)\epsilon(i,1,2,\cdots,\hat i,\cdots
,n-1)\epsilon(I,J)\\
&=&(-1)^{k-1}(-1)^{i-1}=(-1)^{k+i}. \end{eqnarray*}
 Therefore, the
right side of (\ref{eq:399}) is equal to
\begin{eqnarray*}
&&\sum_{I\cap J=\emptyset,(I\cup
J)^c=\{i\}}\epsilon(I,i,J)\epsilon(I',J')\epsilon(I',J')\epsilon(I,J)\epsilon(I,J)\frac{\det
\overline{u}_{J',J}\det \overline{u}_{I',I}}
{\det \overline{u}}\\
&=&(-1)^{k+i}\epsilon(I',J') \sum_{I\cap J=\emptyset,(I\cup
J)^c=\{i\}}\epsilon(I',J')\epsilon(I,J)\frac{\det
\overline{u}_{J',J}\det \overline{u}_{I',I}} {\det \overline{u}}
\end{eqnarray*}

According to Lemma~\ref{lem:det}(iii), the above sum is 0 if
$J'\cap I'\ne\emptyset$. Otherwise, if $J'\cap I'=\emptyset$ so
that $(I'\cup J')^c=\{l\}$,  the right side of (\ref{eq:399})
equals
\[
\frac{(-1)^{k-l}\epsilon(I',J')(-1)^{i+l}}{\det\overline{u}}
\sum_{I\cap J=\emptyset,(I\cup
J)^c=\{i\}}\epsilon(I',J')\epsilon(I,J)\det
\overline{u}_{J',J}\det \overline{u}_{I',I}.
\]

Now by Lemma~\ref{lem:det}(ii), the above sum is the determinant
of the $(l,i)$-minor of the matrix $\overline{u}$, call this
$\overline{M}_{li}$. Thus, for $J'\cap I'=\emptyset$, the right
side of (\ref{eq:399}) is equal to
\begin{eqnarray*}
\epsilon(I',l,J')\left[\frac{(-1)^{i+l}\det\overline{M}_{li}}{\det\overline{u}}\right]
&=&\epsilon(I',l,J')\times \left[ \mbox{ the }(i,l)\mbox{-entry of
the inverse of } \overline{u}\right]\\
&=&\epsilon(I',l,J')u_{li}.\qed
\end{eqnarray*}

\begin{theorem}\label{thm:1}
$H_n^k$ is a homogeneous operator space.
\end{theorem}
\proof\ Let $\alpha$ be a unitary operator on $H_n^k$. To prove
the theorem, it suffices, by \cite[Prop.\ 9.2.1]{Pisier03}, to
show that $\alpha$ is a complete isometry. We shall show that
$\alpha(x)=\lambda vxw$ for suitable unitary matrices $v$ and $w$,
and $\lambda\in\CC$, with $|\lambda|=1$, which will complete the
proof.

Recall that $\psi:C_n\rightarrow H_n^k$ is the isometry defined by
$\psi(e_i)=b_i^{n,k}$.  Let $\psi^{-1}\alpha\psi$ have matrix
$u^{-1}$ on $C_n$ with respect to the basis $e_1,\ldots,e_n$. As
in Lemma~\ref{lem:1.1}, let $u_1,\ldots,u_n$ be the columns of
$u$. We shall show that $\alpha(x)=\lambda vxw$ holds for every
$x\in H_n^k$, where
 $\lambda=\det\overline{u}$,
$ w=\left[\begin{array}{ccc} \wedge_{i\in
I_1}u_i&\cdots&\wedge_{i\in I_q}u_i
\end{array}\right]
$, and
\[
v=\left[\begin{array}{c} (\wedge_{j\in J_1}u_j)^t\\
\vdots \\
 (\wedge_{j\in
J_p}u_j)^t
\end{array}\right].
\]
(The fact that $v$ and $w$ are unitary matrices follows from the
definition of the inner product on $\wedge^r\CC^n$:
$\ip{x_1\wedge\cdots\wedge x_r}{y_1\wedge\cdots\wedge y_r}
=\det[\ip{x_i}{y_j}]$.)

In the first place, $\psi^{-1}\alpha\psi(u_i)=u^{-1}(u_i)=e_i$, so
that $\alpha\psi(u_i)=\psi(e_i)=b_i^{n,k}$. Thus it suffices to
prove
\begin{equation}\label{eq:claim2}
v\psi(u_i)w=b_i^{n,k}/\det\overline{u}.
\end{equation}
Let us first show that
\begin{equation}\label{eq:claim1}
\psi(u_i)=(\det\overline{u})^{-1}\sum_{I\cap J=\emptyset,(I\cup
J)^c=\{i\}}\epsilon(I,i,J) \left(\wedge_{j\in J}\overline
{u}_j\right)\left(\wedge_{i\in I}\overline{u}_i\right)^t.
\end{equation}

By Lemma~\ref{lem:1.1}, the proof of (\ref{eq:claim1}) amounts to
\begin{equation}\label{eq:456}
\left[\left(\wedge_{j\in J}u_j\right)\left(\wedge_{i\in
I}u_i\right)^t\right]_{J',I'}= \det u_{J',J}\det u_{I',I}.
\end{equation}
The left side of (\ref{eq:456}) is given by
$e_{J'}^t\left(\wedge_{j\in J}u_j\right)\left(\wedge_{i\in
I}u_i\right)^te_{I'}$. By Lemma~\ref{lem:det}(i), $\wedge_{j\in
J}u_j=\sum_L\det u_{L,J}e_L$, where $L$ runs over the subsets of
cardinality $n-k$. Hence $e_{J'}^t\left(\wedge_{j\in
J}u_j\right)\\ =\sum_L\det u_{L,J}e_{J'}^te_L=\det u_{J',J}$.
Similarly, $\left(\wedge_{i\in I}u_i\right)^te_{I'}=\det
u_{I',I}$, which proves (\ref{eq:456}) and hence
(\ref{eq:claim1}).

We now use (\ref{eq:claim1}) to prove (\ref{eq:claim2}). Note that
since
$e_{J,I}=e_J\otimes\overline{e}_I=e_J(e_I)^t=\left(\wedge_{j\in
J}e_j\right)\left(\wedge_{i\in I}e_i\right)^t$, we may write
\begin{equation}\label{eq:455}
b_i^{n,k}=\sum_{I\cap J=\emptyset,(I\cup
J)^c=\{i\}}\epsilon(I,i,J)\left(\wedge_{j\in
J}e_j\right)\left(\wedge_{i\in I}e_i\right)^t.
\end{equation}
By (\ref{eq:claim1}) and (\ref{eq:455}), it suffices to prove
\[
v\left(\wedge_{j\in J}\overline {u}_j\right)\left(\wedge_{i\in
I}\overline{u}_i\right)^tw=\left(\wedge_{j\in J}e_j\right)
\left(\wedge_{i\in I}e_i\right)^t.
\]
This is a simple calculation. Suppose for definiteness that
$J=J_r$ and $I=I_s$. Then $ v\left(\wedge_{j\in J}\overline
{u}_j\right)= e_{J_r}$,  and $ \left(\wedge_{i\in
I}\overline{u}_i\right)^tw=e_{I_s}^t$. \qed

\begin{remark}
In {\rm \cite[page 2230]{NeaRus03}}, we defined an operator space
construction
 denoted by $\mbox{Diag}(H_n^{k_1},\ldots,H_n^{k_m})$  which
depended on a choice of orthonormal basis for each of the spaces
$H_n^{k_j}$. Because of the homogeneity of the spaces $H_n^k$
proved in Theorem~\ref{thm:1}, this space is independent of these
choices up to complete isometry and is now seen to be the
intersection $H_n^{k_1}\cap\ldots\cap H_n^{k_m}$ in the sense of
operator space theory {\rm (\cite[page 55]{Pisier03})}.
\end{remark}

%\subsection{This is a numbered second-level section head}
%This is an example of a numbered second-level heading.

%\subsection*{This is an unnumbered second-level section head}
%This is an example of an unnumbered second-level heading.

%\subsubsection{This is a numbered third-level section head}
%This is an example of a numbered third-level heading.

%\subsubsection*{This is an unnumbered third-level section head}
%This is an example of an unnumbered third-level heading.

\section{Anti-symmetric Fock spaces}\label{sec:22}

Let $ C^{n,k}_{h}$ denote the wedge (or creation) operator from
$\wedge^{k-1}\CC^n$ to $\wedge^{k}\CC^n$ given by
$$C^{n,k}_{h}(h_1\wedge\cdots\wedge h_{k-1})=h\wedge
h_1\wedge\cdots\wedge h_{k-1}.
$$

Many properties of these classical operators on the full
anti-symmetric Fock space %$\oplus_{j=0}^{n-1}\wedge^j\CC^n$
 are
given
 in \cite[Exercises 12.4.39-40]{KadRin86}.

\medskip

As in section~\ref{sec:1}, let $e_1,\ldots,e_n$ be the usual
column vector orthonormal basis for $\CC^n$,
 and let $\{e_{I_1},\ldots,e_{I_q}\}$  and $\{e_{J_1},\ldots,e_{J_p}\}$
 be the column vector orthonormal bases for $\CC^q$ and
 $\CC^p$ respectively, and define
the unitary operators $U_j^n$ ($j=k-1$ and $j=n-k$),\ $W_k^n,\
V_k^n$ in the diagram below as follows:
\begin{itemize}
\item $U_{k-1}^n(e_I)=e_{i_1}\wedge\cdots\wedge e_{i_{k-1}}$, where
$I=\{i_1<\cdots<i_{k-1}\}$.
\item $U_{n-k}^n(e_J)=e_{j_1}\wedge\cdots\wedge e_{j_{n-k}}$, where
$J=\{j_1<\cdots<j_{n-k}\}$.
\item $V_k^n(e_{i_1}\wedge\cdots\wedge e_{i_k})=e_{j_1}\wedge\cdots\wedge
e_{j_{n-k}}$, where $\{j_1<\cdots<j_{n-k}\}$ is the complement of
$\{i_1<\cdots<i_k\}$.
\item $W_k^n(e_{j_1}\wedge\cdots\wedge
e_{j_{n-k}})=\epsilon(i,I)\epsilon(I,i,J)e_{j_1}\wedge\cdots\wedge
e_{j_{n-k}}$ for any $i$ and $I$ such that $I\cap J=\emptyset$ and
$(I\cup J)^c=\{i\}$ (which is independent of the choice of $i$ or
$I$).
\end{itemize}

%\bigskip

\[
\begin{array}{ccccc}
\CC^q&\stackrel{b_i^{n,k}}{\longrightarrow} & \CC^p\\
                    &&\\
U_{k-1}^n\downarrow\quad\quad\quad &
&\quad\quad\quad\downarrow U_{n-k}^n\\
&&\\
\wedge^{k-1}\CC^n& & \wedge^{n-k}\CC^n\\
&&\\
 C^{n,k}_{e_i}\downarrow\quad\quad\quad &
&\quad\quad\quad\downarrow W_k^n\\
&&\\
\wedge^k\CC^n&\stackrel{V_k^n}{\longrightarrow}&\wedge^{n-k}\CC^n
\end{array}
\]

\medskip

Note that since $b_i^{n,k}$ is a $p\times q$ matrix, it is viewed
as an operator from $\CC^q$ to $\CC^p$. In the definition of
$W_k^n$, $\epsilon(i,I)$ is the signature of the permutation
$(i,i_1,\ldots,i_{k-1})\mapsto (i_1,\ldots,i,\ldots,i_{k-1})$. To
prove the non-dependence on $i$, suppose $i,i'\not\in J$. Then
\begin{eqnarray*}
\epsilon(I,i,J)&=&\epsilon(i_1,\ldots,i_{k-1},i,j_1,\ldots,j_{n-k})\\
&=&(-1)^{k-1} \epsilon(i,i_1,\ldots,i_{k-1},j_1,\ldots,j_{n-k})\\
&=&(-1)^{k-1}\epsilon(i,I)
\epsilon(i_1,\ldots,i,\ldots,i_{k-1},j_1,\ldots,j_{n-k}),
\end{eqnarray*}
where $i_1<\cdots <i<\cdots i_{k-1}$. Similarly,
\[
\epsilon(I',i',J)=(-1)^{k-1}\epsilon(i',I')
\epsilon(i'_1,\ldots,i',\ldots,i'_{k-1},j_1,\ldots,j_{n-k}),
\]
where $i'_1<\cdots <i'<\cdots i'_{k-1}$. Hence,
$\epsilon(i,I)\epsilon(I,i,J)=\epsilon(i',I')\epsilon(I',i',J)$.

\medskip

It is now a simple matter to check that
$W_k^nU_{n-k}^nb_i^{n,k}=V_k^nC_{e_i}^{n,k}U_{k-1}^n $. Indeed,
for any $I'$, $
W_kU_{n-k}^nb_i^{n,k}(e_{I'})=W_kU_{n-k}^n\epsilon(I',i,J)(e_J)=
W_k(\epsilon(I',i,J)e_{j_1}\wedge\cdots\wedge e_{j_{n-k}}=
\epsilon(i,I')e_{j_1}\wedge\cdots\wedge e_{j_{n-k}}, $
 and
$ V_k^nC_{e_i}^{n,k} U_{k-1}^n (e_{I'})=V_k^n(e_i\wedge
e_{i'_1}\wedge\cdots\wedge e_{i'_{k-1}})
=V_k^n(\epsilon(i,I')e_{i'_1}\wedge\cdots\wedge
e_i\wedge\cdots\wedge e_{i'_{k-1}})=
\epsilon(i,I')e_{j_1}\wedge\cdots\wedge e_{j_{n-k}}$.

Hence, letting ${\mathcal C}^{n,k}$ denote the space
$\mbox{sp}\{C_{e_i}^{n,k}\}$ yields the following lemma.
\begin{lemma}\label{lem:2.1}
$H_n^k$ is completely isometric to ${\mathcal C}^{n,k}$.
\end{lemma}

By \cite[Theorem 2, Corollary 2.8]{NeaRus03}, every atomic
contractively complemented subspace $X$ of a \csa\ is
isometrically completely contractive to a direct sum of Cartan
factors of types 1 to 4 and some of the spaces $H_n^k$.  The
following theorem gives more detailed information on what can be
said up to complete isometry in the case of an Hilbertian $X$.

Recall that a linear map of one operator space into another is
said to be a {\it complete semi-isometry} if it is isometric and
completely contractive.

\begin{theorem}\label{thm:2}
Let $X$ be the range of a contractive projection $P$ on a \csa\
$A$, and suppose that $X$ is isometric to a Hilbert space. Then
there exist projections $p,q\in A^{**}$ such that
\begin{description}
\item[(a)] $X=\{pxq+(1-p)x(1-q):x\in X\}$;
\item[(b)] The map ${\mathcal E}_0x=pxq$ is a complete semi-isometry of $X$ onto $pXq$;
\item[(c)] If $X$ is finite-dimensional, then $pXq$ is completely isometric to an intersection
of the spaces ${\mathcal C}^{n,k}$. If $X$ is
infinite-dimensional, then $pXq$ is completely semi-isometric to
either row or column Hilbert space;
\item[(d)] Both $X$ and $pXq$ are completely isometric to the range of a contractive
projection on $B(K)$ for an appropriate Hilbert space $K$.
\item[(e)] $P^{**}(pxq)=x$ and $P^{**}((1-p)x(1-q))=0$, for $x\in P(A)$. Hence
 $P^{**}:pXq\rightarrow
X$ is the inverse of ${\mathcal E}_0$ and $(1-p)X(1-q)\subset \ker
P^{**}$.
\end{description}
\end{theorem}
\proof\ \begin{description}
\item[(a)] Since $P(A)$ is reflexive, $X=P(A)=P^{**}(A^{**})=pXq+(1-p)X(1-q)$, the last
equality following from \cite[Prop. 4]{FriRus85};
\item[(b)] By \cite[Lemma 2.2]{NeaRus03};
\item[(c)] By Lemma~\ref{lem:2.1} and \cite[Prop. 2.6]{NeaRus03};
\item[(d)] For $X$ this follows directly from \cite[Corollary 2.8]{NeaRus03}, and for
 $pXq$ it follows directly from
\cite[Theorem 3(b) and Corollary 7.3]{NeaRus03}. \footnote{The
authors wish to take this opportunity to point out the following
correction to
  \cite[Lemma 7.2 and Corollary 7.3]{NeaRus03}. The term $
{n-1 \choose k-1}^{1/2}$ should be replaced by ${n-1 \choose k-1}$
in the statements of Lemma 7.2 and Corollary 7.3, and in the proof
of Lemma 7.2. Accordingly, in the proof of Corollary 7.3,
$m^{1/2}$ should be replaced by $m$.}
\item[(e)] For $x\in P(A)$, $P^{**}(pxq)\in P^{**}(A^{**})=P(A)$, say
$P^{**}(pxq)=y\in P(A)$ and by \cite[Prop. 4]{FriRus85},
$y=pyq+(1-p)y(1-q)$. Thus $p(P^{**}(pxq))q=pyq$. Now it follows
from \cite[Prop. 1]{FriRus85}, that $p(P^{**}z)q=z$ for all $z\in
pP^{**}(A^{**})q$. With $z=pxq$, we have $pxq=pyq$ and by (b)
$x=y$, proving that $P^{**}(pxq)=x$. Moreover,
$P^{**}((1-p)x(1-q))=P^{**}(x-pxq)=x-P^{**}(pxq)=x-y=0$. \qed
\end{description}

\begin{remark}

For any subset $S\subset\{1,\ldots,n\}$, $\cap_{k\in S}{\mathcal
C}^{n,k}$ is exactly the space of creation operators on the direct
sum $\oplus_{k\in S}\wedge^{k-1}\CC^n$. So Theorem~\ref{thm:2}
says that all finite-dimensional contractively complemented
Hilbertian operator spaces are essentially a space of creation
operators in the anti-symmetric Fock space. In particular, if
$S=\{1,\ldots,n\}$, this is the space $\Phi_n$ discussed in {\rm
\cite[section 9.3]{Pisier03}}.
\end{remark}

The following two properties of the wedge operators
 follow easily from \cite[Exercises 12.4.39-40]{KadRin86}. They will be used, together
 with Lemma~\ref{lem:2.2}, in section~\ref{sec:2}.

\begin{eqnarray}\nonumber
C_h^{n,k*}C_h^{n,k}(h_1\wedge\cdots\wedge
h_{k-1})&=&\ip{h}{h}h_1\wedge\cdots\wedge h_{k-1}
-\ip{h_1}{h}h\wedge h_2\wedge\cdots\wedge h_{k-1}\\
\label{eq:1} &&+\cdots \pm\ip{h_{k-1}}{h}h\wedge
h_1\wedge\cdots\wedge h_{k-2}.
\end{eqnarray}
and
\[
C_h^{n,k}C_h^{n,k*}(h_1\wedge\cdots\wedge
h_k)=\sum_{j=1}^k\ip{h_j}{h} h_1\wedge\cdots\wedge
h'_j\wedge\cdots \wedge h_k\quad (h'_j=h).
\]  In particular,
$C_h^{n,1}C_h^{n,1*}=h\otimes\overline{h}$, for $h\in \CC^n$.

\medskip

\begin{lemma}\label{lem:2.2}
\begin{description}
\item[(a)] $\mbox{tr}(C_h^{n,k*}C_h^{n,k})={{n-1}\choose{k-1}}\|h\|^2$.
In particular, $C_h^{n,1*}C_h^{n,1}=\|h\|^2$.
\item[(b)] Let the (repeated) eigenvalues of $\sum_{i=1}^mC_{h_i}^{n,1}
C_{h_i}^{n,1*}$ be $\lambda_1\ge\cdots\ge\lambda_n$. Then the
eigenvalues of $\sum_{i=1}^mC_{h_i}^{n,k} C_{h_i}^{n,k*}$ are
precisely the sums of $k$ eigenvalues of
$\sum_{i=1}^mC_{h_i}^{n,1} C_{h_i}^{n,1*}$.
\end{description}
\end{lemma}
\proof\ In the first place, we have
\begin{eqnarray*}
\mbox{tr}(C_h^{n,k}C_h^{n,k*})&=&\mbox{tr}(C_h^{n,k*}C_h^{n,k})\\
&=&\sum_{i_1,\ldots,i_{k-1}}\ip{C_h^{n,k*}C_h^{n,k}(e_{i_1}\wedge\cdots
\wedge e_{i_{k-1}})}{e_{i_1}\wedge\cdots \wedge e_{i_{k-1}}}\\
&=&\sum_{i_1,\ldots,i_{k-1}}\ip{h\wedge e_{i_1}\wedge\cdots
\wedge e_{i_{k-1}}}{h\wedge e_{i_1}\wedge\cdots \wedge e_{i_{k-1}}}\\
&=&\sum_{i_1,\ldots,i_{k-1}}\det\left[\begin{array}{cccc}
\ip{h}{h}&\ip{h}{e_{i_1}}&\cdots&\ip{h}{e_{i_{k-1}}}\\
\ip{e_{i_1}}{h}&1&\cdots&0\\
\cdots&\cdots&\cdots&\cdots\\
\ip{e_{i_{k-1}}}{h}&0&\cdots&1
\end{array}\right]\\
&=&\sum_{i_1,\ldots,i_{k-1}}\left(\|h\|^2-\sum_{j=1}^{k-1}|\ip{h}{e_{i_j}}|^2\right)\\
&=&{n\choose k-1}\|h\|^2-\sum_{i_1,\ldots,i_{k-1}}\sum_{j=1}^{k-1}|\ip{h}{e_{i_j}}|^2\\
&=&{n\choose k-1}\|h\|^2-\sum_{l=1}^n{n-1\choose k-2}|\ip{h}{e_l}|^2\\
&=& \left( {n\choose k-1}-{n-1\choose
k-2}\right)\|h\|^2={n-1\choose k-1}\|h\|^2.
\end{eqnarray*}

To prove the second statement, let $\xi_1,\ldots,\xi_n$ be an
orthonormal basis of $\CC^n$ consisting of eigenvectors of
$\sum_{i=1}^mC_{h_i}^{n,1} C_{h_i}^{n,1*}$ so that
$\sum_i\ip{\xi_k}{h_i}h_i=\lambda_k\xi_k$. Then
\begin{eqnarray*}
\sum_{i=1}^mC_{h_i}^{n,k}
C_{h_i}^{n,k*}(\xi_{i_1}\wedge\cdots\wedge\xi_{i_k})&=&\sum_i\sum_{j=1}^k\ip{\xi_{i_j}}{h_i}
\xi_{i_1}\wedge\cdots\wedge h_i\wedge\cdots\wedge\xi_{i_k}\\
&=&\sum_j\xi_{i_1}\wedge\cdots\wedge\left[\sum_i\ip{\xi_{i_j}}{h_i}h_i\right]\wedge\cdots\xi_{i_k}\\
&=&\left[\sum_j\lambda_{i_j}\right]\xi_{i_1}\wedge\cdots\wedge\xi_{i_k}.
\end{eqnarray*}

Conversely, if
$\xi=\sum\alpha_{i_1,\ldots,i_k}\xi_{i_1}\wedge\cdots\wedge\xi_{i_k}$
is an eigenvector of $\sum_{i=1}^mC_{h_i}^{n,k} C_{h_i}^{n,k*}$,
with eigenvalue $\lambda$, then
\begin{eqnarray*}\lefteqn{\lambda
\sum\alpha_{i_1,\ldots,i_k}\xi_{i_1}\wedge\cdots\wedge\xi_{i_k}=}\\
&=& \sum_{i=1}^mC_{h_i}^{n,k}
C_{h_i}^{n,k*}\sum\alpha_{i_1,\ldots,i_k}\xi_{i_1}\wedge\cdots\wedge\xi_{i_k}\\
&=&\sum\alpha_{i_1,\ldots,i_k}\sum_i\sum_j\ip{\xi_{i_j}}{h_i}\xi_{i_1}\wedge\cdots\wedge
h_i\wedge\cdots \wedge\xi_{i_k}\\
&=&\sum\alpha_{i_1,\ldots,i_k}\sum_j
\xi_{i_1}\wedge\cdots\wedge\left[\sum_{i=1}^mC_{h_i}^{n,1}
C_{h_i}^{n,1*}\xi_{i_j}\right]\wedge\cdots \wedge\xi_{i_k}\\
&=&\sum\alpha_{i_1,\ldots,i_k}\left[\sum_j\lambda_{i_j}\right]\xi_{i_1}
\wedge\cdots\wedge\xi_{i_k}.
\end{eqnarray*}
From this, the second statement follows.\qed
\section{Completely bounded Banach-Mazur distance}\label{sec:2}

Recall that the completely bounded Banach-Mazur distance between
two (completely isomorphic) operator spaces $E,F$ is defined by
\[
\mbox{d}_{\mbox{cb}}(E,F)=\inf\{\|u\|_{\mbox{cb}}
\|u^{-1}\|_{\mbox{cb}}:u:E\rightarrow F\mbox{ complete isomorphism
}\}.
\]

We shall explicitly compute $\mbox{d}_{\mbox{cb}}(H_n^k,H_n^1)$.
 By Lemma~\ref{lem:2.1}, we can
identify $H_n^k$ with ${\mathcal C}^{n,k}$. For a fixed $k$, let
$\psi:H_n^1\rightarrow H_n^k$ be the isometry given by
$\psi(C_{e_i}^{n,1})=C_{e_i}^{n,k}$. Recall from \cite[Prop. 4]
{Mathes94} that, since $\psi$ is a mapping from the column Hilbert
space $H_n^1$, $\|\psi\|_{\mbox{cb}}=\|\psi\|_{\mbox{row-cb}}$,
where
$\|\psi\|_{\mbox{row-cb}}:=\sup\{\|(\psi(h_1),\ldots,\psi(h_m))\|\}$
where the supremum is extended over all $m\ge1$ and all row
vectors with $\|(h_1,\ldots,h_m)\|\le 1$. The norm
$\|\psi\|_{\mbox{col-cb}}$ is defined analogously and by
 \cite[Prop. 2]{Mathes94},
$\|\psi^{-1}\|_{\mbox{cb}}=\|\psi^{-1}\|_{\mbox{col-cb}}$.

\begin{lemma}\label{lem:2}
$\|\psi\|_{\mbox{row-cb}}=\sqrt{k}$ and
$\|\psi^{-1}\|_{\mbox{col-cb}}=\sqrt{\frac{n}{n-k+1}}$.
\end{lemma}
\proof\ Let $A=(C_{h_1}^{n,1},\ldots,C_{h_m}^{n,1})$ and
$B=(C_{h_1}^{n,k},\ldots,C_{h_m}^{n,k})$. We show first that
$\|B\|\le\sqrt k\, \|A\|$. We have $AA^*=\sum_{i=1}^mC_{h_i}^{n,1}
C_{h_i}^{n,1*}=\sum_{i=1}^mh_i\otimes\overline{h_i}$ and
$BB^*=\sum_{i=1}^mC_{h_i}^{n,k} C_{h_i}^{n,k*}$.

Let the (repeated) eigenvalues of $AA^*$ be
$\lambda_1\ge\cdots\ge\lambda_n$. By Lemma~\ref{lem:2.2}(b), we
have
\[
\|B\|^2=\|BB^*\|=\lambda_1+\cdots+\lambda_k\le
k\lambda_1=k\|AA^*\|=k\|A\|^2.
\]
Taking $m=n$ and $h_j=e_j$, we have $AA^*=I$ and $BB^*=kI$,
proving the first statement of the Lemma.

Let $D=(C_{h_1}^{n,1},\ldots,C_{h_m}^{n,1})^t$ and
$C=(C_{h_1}^{n,k},\ldots,C_{h_m}^{n,k})^t$. We show next that
$\|D\|\le\sqrt{\frac{n}{n-k+1}}\|C\|$. By Lemma~\ref{lem:2.2}(a),
we have $D^*D=\sum_{i=1}^mC_{h_i}^{n,1*}
C_{h_i}^{n,1}=\sum_{i=1}^m\|h_i\|^2$ and
$C^*C=\sum_{i=1}^mC_{h_i}^{n,k*} C_{h_i}^{n,k}$. Since $C^*C$ is a
square matrix of size $n\choose k-1$, again by
Lemma~\ref{lem:2.2}(a),
\[
\|C^*C\|{n\choose k-1}\ge \mbox{tr}(C^*C)={n-1\choose k-1}
\sum_{i=1}^m\|h_i\|^2.
\]
Therefore,
\[
\frac{\|D\|^2}{\|C\|^2}=\frac{\sum_{i=1}^m\|h_i\|^2}{\|C^*C\|}
\le\frac{\sum_{i=1}^m\|h_i\|^2} {\frac{{n-1\choose k-1}}{{n\choose
k-1}}\sum_{i=1}^m\|h_i\|^2}=\frac{n}{n-k+1}.
\]
Taking $m=n$ and $h_i=e_i$, we have $D^*D=n$. By (\ref{eq:1}),
$C_{e_i}^{n,k*}C_{e_i}^{n,k}(e_{i_1}\wedge\cdots\wedge
e_{i_{k-1}})=0$ if $i\in \{i_1,\ldots,i_{k-1}\}$ and equal to
$e_{i_1}\wedge\cdots\wedge e_{i_{k-1}}$ otherwise. Hence
$C^*C=(n-k+1)I$, proving the second statement.\qed

\medskip

\begin{theorem}\label{thm:3}
$\mbox{d}_{\mbox{cb}}(H_n^k,H_n^1)=\sqrt{\frac{kn}{n-k+1}}$, for
$1\le k\le n$.
\end{theorem}
\proof\ By \cite[Theorem 3.1]{Zhang97}, and the first paragraph of
its proof,
$\mbox{d}_{\mbox{cb}}(H_n^k,H_n^1)=\|\psi\|_{\mbox{cb}}\|\psi^{-1}\|_{\mbox{cb}}$.
Now apply Lemma~\ref{lem:2} and the remarks just preceding it.
\qed

\medskip

Not surprisingly, we obtain the result published first by Mathes
(\cite[Prop. 7]{Mathes94}, \cite [p. 21]{Pisier03}).

\begin{corollary}
$\mbox{d}_{\mbox{cb}}(R_n,C_n)=n$.
\end{corollary}

Symmetry considerations in Theorem~\ref{thm:3} suggest
$\mbox{d}_{\mbox{cb}}(H_n^k,H_n^n)=
\mbox{d}_{\mbox{cb}}(H_n^{n-k+1},H_n^1)$, which can be proved by
exactly the same methods. Hence we obtain the following, which is
the answer to Problem 1 in \cite{NeaRus03}.

\begin{corollary}
$\mbox{d}_{\mbox{cb}}(H_n^k,H_n^n)=\sqrt{\frac{(n-k+1)n}{k}}$, for
$1\le k\le n$.
\end{corollary}

\begin{remark}
It is curious to note that by Theorem~\ref{thm:3} and its two
corollaries, $\mbox{d}_{\mbox{cb}}(H_n^k,H_n^n)
\mbox{d}_{\mbox{cb}}(H_n^k,H_n^1)=
\mbox{d}_{\mbox{cb}}(H_n^n,H_n^1)$, so that in the metric space of
all operator spaces of dimension $n$ (\cite[page 335]{Pisier03}),
the three points $H_n^k, H_n^1, H_n^n$ form a degenerate triangle.
\end{remark}

\medskip

Since $H_n^k$ is distinct from row or column Hilbert space, new
ideas will be needed to solve the following problem.
\begin{problem}
Find $\mbox{d}_{\mbox{cb}}(H_n^{k_1},H_n^{k_2})$ for
$1<k_1<k_2<n$.
\end{problem}

We have already mentioned in the introduction that $C_n^*=R_n$ and
$R_n^*=C_n$ in the category of operator spaces.  Hence $H_n^1$ and
$H_n^n$ are operator space duals of each other. The following
problem is therefore of interest and its solution would certainly
lead to insight into Pisier's question on the operator space dual
of $\Phi_n$, \cite[page 175]{Pisier03}.

\begin{problem}
Find the operator space dual of $H_n^k$.
\end{problem}

Note that Theorem~\ref{thm:2} does not say anything about the
infinite-dimensional case up to complete isometry, although in
this case, the space is completely semi-isometric to $R$ or $C$.

\begin{problem}\label{prob:3}
Are all infinite-dimensional Hilbertian contractively complemented
operator spaces completely isometric to a space of creation
operators on a subspace of the anti-symmetric Fock space?
\end{problem}

\begin{problem}
What is the completely bounded Banach-Mazur distance between two
infinite-dimensional Hilbertian contractively complemented
operator spaces?
\end{problem}

\bibliographystyle{amsplain}

\begin{thebibliography}{10}

\bibitem{AraFri78} Jonathan Arazy and Yaakov Friedman, {\it Contractive projections
in $C_1$ and $C_\infty$}, Mem. Amer. Math. Soc. 13 (1978), no 200.

\bibitem{BleLeM04} David P. Blecher and Christian Le Merdy, Operator algebras and their
modules---an operator space approach, Clarendon Press, Oxford
(forthcoming)

\bibitem{EffRau00} Edward Effros and Zhong-Jin Ruan, Operator Spaces, Oxford University
Press, 2000.

\bibitem{FriRus85} Yaakov Friedman and Bernard Russo, {\it Solution of
the contractive projection problem}, J. Funct. Anal. 60 (1985),
56--79.

\bibitem{KadRin86} Richard V. Kadison and John R. Ringrose,
Fundamentals of the Theory of Operator Algebras, Volume II,
Academic Press 1986

\bibitem{Mathes94} Ben Mathes, {\it Characterizations of row and column Hilbert space}, J. London
Math. Soc. (2) {\bf 50} (1994) 199--208.

\bibitem{NeaRus03} Matthew Neal and Bernard Russo, {\it Contractive projections and operator
spaces}, Trans. Amer. Math. Soc. {\bf 355} (2003), 2223--2362.

\bibitem{Paulson03} Vern Paulson, Completely bounded maps and operator algebras,
Cambridge Studies in Advanced Mathematics 78, Cambridge University
Press, Cambridge, 2002.

\bibitem{Pisier03} Gilles Pisier, Introduction to Operator Space Theory, Cambridge University Press
2003.

\bibitem{Robertson91} A. Guyan Robertson, {\it Injective matricial Hilbert spaces}, Math. Proc.
Cambridge Philos. Soc. 110 (1991), 183--190.

\bibitem{Yokonuma92} Takeo Yokonuma, Tensor spaces and exterior algebra, Translations of Mathematical
Monographs, Volume 108, American Mathematical Society 1992.

\bibitem{Zhang97} Chun Zhang, {\it Completely bounded Banach-Mazur distance}, Proc.
Edinburgh Math. Soc. {\bf 40} (1997) 247--260.

\end{thebibliography}

\end{document}